\documentclass[number,dvips]{arxbj}

\aid{0}
\volume{16}
\issue{2}
\pubyear{2010}
\firstpage{561}
\lastpage{584}
\doi{10.3150/09-BEJ209}

\makeatletter
\newtheorem{corollary}{Corollary}[section]
\newtheorem{theorem}{Theorem}[section]
\newtheorem{proposition}{Proposition}[section]
\newremark{example}{Example}[section]
\newtheorem{lemma}{Lemma}[section]

\newremark{remark}{Remark}[section]

\makeatother

\begin{document}
\begin{frontmatter}

\title{Uniform error bounds for a continuous approximation of
non-negative random~variables}
\runtitle{Uniform error bounds}

\begin{aug}
\author{\fnms{Carmen} \snm{Sang\"{u}esa}\ead[label=e1]{csangues@unizar.es}}
\runauthor{C. Sang\"{u}esa}
\pdfauthor{Carmen Sanguesa}
\address{Departamento de M\'{e}todos Estad\'{\i}sticos,
Facultad de Ciencias, Universidad de Zaragoza,
Pedro Cerbuna, 12, 50009 Zaragoza, Spain. \printead{e1}}
\end{aug}

\received{\smonth{11} \syear{2008}}

%
\begin{abstract}
In this work, we deal with approximations for distribution functions of
non-negative
random variables. More specifically, we construct continuous
approximants using an acceleration technique over a well-know
inversion formula for Laplace transforms. We give uniform error
bounds using a representation of these approximations in terms of
gamma-type operators. We apply our results to certain mixtures of
Erlang distributions which contain the class of continuous
phase-type distributions.
\end{abstract}

%
\begin{keyword}
\kwd{gamma distribution}
\kwd{Laplace transform}
\kwd{phase-type distribution}
\kwd{uniform distance}
\end{keyword}

\end{frontmatter}

\section{Introduction}\label{sec1}

Frequent operations in probability such as convolution or random
summation of random variables produce probability distributions
which are difficult to evaluate in an explicit way. In these cases,
one needs to use numerical evaluation methods. For instance, one can
use numerical inversion of the Laplace or Fourier transform of the
distribution at hand (see \cite{abwhan} for the general use of
Laplace--Stieltjes transforms in applied probability or \cite
{emgrso,grheco} for the method of Fast Fourier
Transform in the context of risk theory). Another approach is the use
of recursive evaluation methods, of special interest for random sums
(see \cite{grheco,surecu}, for instance). Some of the methods mentioned
above require a previous discretization step
to be applied to the initial random variables when these are
continuous. The
usual way to do so is by means of rounding methods. However, it is
not always possible to evaluate the distribution of the rounded random
variable in an
explicit way and it is not always clear when using these methods how
the rounding error propagates when one takes successive
convolutions. In these cases, it seems worthwhile to consider
alternative discretization methods. For instance, when dealing with
non-negative random variables, the following method (\cite{feanin}, page
233) has been proposed in the literature. Let $X$ be a
random variable taking values in $[0,\infty)$
with distribution function $F$. Denote
by $\phi_{X}(\cdot)$ the Laplace--Stieltjes (LS) transform of $X$, that
is,
\[
\phi_{X}(t):=E\mathrm{e}^{-tX}=\int_{[0,\infty)}\mathrm{e}^{-t u}\,
\mathrm{d}F(u),\qquad t>0.
\]
For each $t>0$, we define a random
variable $X^{\bullet t}$ taking values on $k/t, k\in\mathbb{N}$, and
such that
%
\begin{equation}\label{discre}
P(X^{\bullet t}=k/t)=\frac{(-t)^{k}}{k!}\phi
_{X}^{(k)}(t),\qquad k\in\mathbb{N},
\end{equation}
where $\phi^{(k)}_{X}$ denotes the
$k$th derivative ($\phi^{(0)}_{X}\equiv\phi_{X}$).

Thus, if we denote by $L_{t}^{*}F$ the distribution function of
$X^{\bullet t}$, we have
%
\begin{equation}\label{dfdfdi}
L_{t}^{*}F(x):=P(X^{\bullet t}\leq
x)=\sum_{k=0}^{[tx]}\frac{(-t)^{k}}{k!}\phi^{(k)}_{X}(t),\qquad
x\geq0,
\end{equation}
where $[x]$ indicates the largest
integer less than or equal to $x$. The use of this method allows one to
obtain the
probability mass function in an explicit way in some situations in which
rounding methods could not (see, for instance, \cite{adcaon} for
gamma distributions). Moreover, this method allows for an easy
representation of $L_{t}^{*}F$ in terms of $F$, which makes
possible the study of rates of convergence in the approximation
(\cite{adcaon,adcaap}). In \cite{adcaon}, the problem was studied in
a general setting, whereas in \cite{adcaap}, a detailed analysis was
carried out for the case of gamma distributions, that is, distributions whose
density function is given by
%
\begin{equation}\label{gamden}
f_{a,p}(x):=\frac{a^{p} x^{p-1}\mathrm{e}^{-a
x}}{\Gamma(p)},\qquad x>0.
\end{equation}
Also, in
\cite{saerro}, error bounds for random sums of mixtures of
gamma distributions were obtained, uniformly controlled on the
parameters of the
random summation index. In all of these papers, the measure of distance
considered was the Kolmogorov (or sup-norm) distance. More
specifically, for a given real-valued function $f$ defined on
$[0,\infty)$,
we denote by $\|f\|$ the sup-norm, that is,
\[
\|f\|:=\sup_{x\geq0} |f(x)|.
\]
It was shown in \cite{adcaap} that for gamma distributions with
shape parameter $p\geq1$, we have that $\|L_{t}^{*}F-F\|$ is of
order $1/t$, the length of the discretization interval. Note that
$\|L_{t}^{*}F-F\|$ is the Kolmogorov distance between $X$ and
$X^{\bullet t}$, as both are non-negative random variables.

The aim
of this paper is twofold. First, we will consider a
continuous modification of (\ref{dfdfdi}) and give conditions under
which this continuous
modification has rate of convergence of $1/t^2$ instead of $1/t$
(see Sections \ref{sethea} and \ref{seerro}).
In Section \ref{seapp1}, we will consider
the case of gamma distributions to see that the error bounds are
also uniform on the shape parameter. Finally, in Section
\ref{seapp2}, we will consider the application of the results in
Section \ref{seapp1} to the class of mixtures of Erlang
distributions, recently studied in \cite{wiwoon}. This class
contains many of the distributions used in applied probability (in
particular, phase-type distributions) and is closed under important
operations such as mixtures, convolution and compounding.

\section{The approximation procedure}\label{sethea}

The representation of $L_{t}^{*}F$ in (\ref{dfdfdi}) in terms of a
Gamma process (see \cite{adcaon}) will play an important role in our
proofs. We recall this representation. Let $(S(u) ,u \geq0 ) $ be a gamma
process, in which $S(0)=0$ and such that for $u>0$, each $S(u)$ has a gamma
density with parameters $a=1$ and $p=u$, as given in (\ref{gamden}).
Let $g$ be a function defined on $[0,\infty)$. We consider the
gamma-type operator $L_{t}$ given by
%
\begin{equation}\label{operbue}
L_{t}g(x):= Eg\biggl(\frac{S(tx)}{t}\biggr),\qquad x\geq0, t>0,
\end{equation}
provided that this operator is well defined, that is,
$L_{t}|g|(x)<\infty, x\geq0, t>0$. Then, for $F$
continuous on $(0,\infty)$, $L_{t}^{*}F$ in
(\ref{dfdfdi}) can be written as (see \cite{adcaon}, page 228)
%
\begin{equation}\label{rlftlt}
L_{t}^{*}F(x)=L_{t}F \biggl(\frac{[tx]+1}{t} \biggr)=EF \biggl(\frac
{S([tx]+1)}{t} \biggr),\qquad x\geq0, t>0.
\end{equation}
It can be seen that
the rates of convergence of $L_{t}g$ to $g$ are, at most, of order
$1/t$ (see (\ref{taylor}) below). Our aim now is to get faster
rates of convergence. To this end, we will consider the following
operator, built using a classical acceleration technique (\textit
{\textit{Richardson's extrapolation}} -- see, for instance, \cite
{emgrso,grheco}):
%
\begin{equation}\label{acgamo}
L^{[2]}_{t}g(x):=2L_{2t}g(x)-L_{t}g(x)=2Eg \biggl(\frac
{S(2tx)}{2t} \biggr)-Eg \biggl(\frac{S(tx)}{t} \biggr),\qquad x\geq0.
\end{equation}
We will obtain a rate of uniform convergence from
$L_{t}^{[2]}g$ to $g$, of order $1/t^{2}$, on the following class of
functions:
%
\begin{equation}\label{dfclas}
\mathcal{D}:=\{g\in C^{4}([0,\infty))\dvt
\|x^{2}g^{iv}(x)\|<\infty\}.
\end{equation}
The
problem with $L^{[2]}_{t}g$ is that when $tx$ is not a natural
number, $L_{t}g(x)$ is given in terms of Weyl fractional derivatives
of the Laplace transform (see \cite{adsadi}, page 92) and, in general,
we are not able to compute them in an explicit way. However, if we
modify $L^{[2]}_{t}g$ using linear interpolation, that is,
%
\begin{equation}\label{a2gamo}
M^{[2]}_{t}g(x):=(tx-[tx]) \biggl(L^{[2]}_{t}g \biggl(\frac
{[tx]+1}{t} \biggr) \biggr)+([tx]+1-tx) \biggl(L^{[2]}_{t}g \biggl(\frac{[tx]
}{t} \biggr) \biggr),
\end{equation}
then we observe that the
order of convergence of $M_{t}^{[2]}g$ to $g$ is also $1/t^{2}$ on
the following class of functions:
%
\begin{equation}\label{dfcla2}
\mathcal{D}_{1}:=\{g\in C^{4}([0,\infty))\dvt \|g''(x)\|\leq\infty
\mbox{ and }
\|x^{2}g^{iv}(x)\|<\infty\}.
\end{equation}

Moreover, the advantage of using $M^{[2]}_{t}g$ instead of
$L^{[2]}_{t}g$ to approximate $g$ is the computability. In the
following result, we note that the last approximation applied to a
distribution function $F$ is related to $L_{t}^{*}F$, as defined in
(\ref{dfdfdi}). From now on, $\mathbb{N}^{*}$ will denote the set
$\mathbb{N}
\setminus\{0\}$.

\begin{proposition}\label{dftilf}
Let $X$ be a non-negative random variable with
Laplace transform $\phi_{X}$. Let $L_{t}^{*}F, t>0$, be as defined
in \textup{(\ref{dfdfdi})} and let $M_{t}^{[2]}F$ be as defined in
\textup{(\ref{a2gamo})}. We have
\begin{equation}\label{dftile}
M_{t}^{[2]}F \biggl(\frac{k}{t} \biggr)=
\cases{
F(0), & \quad if $k=0$, \cr
\displaystyle
2L_{2t}^{*}F \biggl(\frac{2k-1}{2 t} \biggr)-L_{t}^{*}F \biggl(\frac{k-1}{t} \biggr), &
\quad if $k\in\mathbb{N}^{*}$, }
\end{equation}
and
%
\begin{equation}\label{dftilt}
M^{[2]}_{t}F(x)=(tx-[tx])M^{[2]}_{t}F \biggl(\frac
{[tx]+1}{t} \biggr)+([tx]+1-tx)M^{[2]}_{t}F \biggl(\frac{[tx]
}{t} \biggr).
\end{equation}
\end{proposition}

\begin{pf}
Let $t>0$ be fixed. First, observe that by (\ref{a2gamo}), we can write
%
\begin{equation}\label{dftil1}
M^{[2]}_{t}F \biggl(\frac{k}{t} \biggr)=L^{[2]}_{t}F \biggl(\frac{k}{t}
\biggr),\qquad k\in\mathbb{N}.
\end{equation}
Now, using (\ref{acgamo}) and (\ref{operbue}), we have
$M_{t}^{[2]}F(0)=
L^{[2]}_{t}F(0)=F(0)$, which shows (\ref{dftile}) for $k=0$. Finally,
using (\ref{acgamo}),
(\ref{operbue}) and (\ref{rlftlt}), we have, for $k\in\mathbb{N}^{*},$
%
\begin{equation}\label{valuek}
L^{[2]}_{t}F \biggl(\frac{k}{t} \biggr)=2EF \biggl( \frac{S(2k)}{
2t} \biggr) -EF \biggl( \frac{S(k)}{t}
\biggr)=2L_{2t}^{*}F \biggl( \frac{2k-1}{2
t} \biggr) -L_{t}^{*}F \biggl( \frac
{k-1}{t} \biggr).
\end{equation}
Thus, (\ref{dftil1}) and (\ref{valuek}) show (\ref{dftile}) for
$k\in
\mathbb{N}^{*}$. Note that
(\ref{dftilt}) is obvious by (\ref{a2gamo}) and (\ref{dftil1}).
This
completes the proof of Proposition \ref{dftilf}.
\end{pf}

In the following example, we illustrate the use of the previous
approximant in the context of random sums, defined in the following
way. Let $(X_{i})_{i\in\mathbb{N}^{*}} $ be a sequence of
independent, identically distributed
non-negative random variables. Let $M$ be a random variable
concentrated on the non-negative integers, independent of
$(X_{i})_{i\in\mathbb{N}^{*}} $.
Consider the random variable
%
\begin{equation}\label{ransum}
\sum_{i=1}^{M}X_{i},
\end{equation}
with the convention that the empty sum is 0.

\begin{example}\label{ex1}
As pointed out in the \hyperref[sec1]{Introduction}, an
explicit expression for the distribution of (\ref{ransum}) is usually
not possible.
Our aim is to consider an example in which this distribution can be
evaluated explicitly and to compare our approximation method with some
others considered in the literature. To this end, we consider that $M$
follows a geometric distribution of parameter $p$, that is,
$P(M=k)=(1-p)^{k}p, k\in\mathbb{N}$ and $(X_{i})_{i\in\mathbb
{N}^{*}} $ are exponentially
distributed (with mean $1$, for the sake of simplicity). In this case,
it is well known (use LS transforms, for instance) that (\ref{ransum})
has the same distribution as a mixture of the degenerate distribution
at 0 (with probability $p$) and an exponential distribution, that is,
%
\begin{equation}\label{testfi}
F(x):=P\Biggl(\sum_{i=1}^{M}X_{i}\leq
x\Biggr)=p+(1-p)(1-\mathrm{e}^{-px})=1-(1-p)\mathrm{e}^{-px},\qquad x\geq0.
\end{equation}
When an explicit expression is not possible, the usual approximate
evaluation method is by discretizing the summands in (\ref{ransum}) and
then using recursive methods found in the literature for discrete
random sums. By considering (\ref{discre}) as a first discretization
method, we have (see \cite{adcaap}, page~391)
%
\begin{equation}\label{gapom2}
P \biggl(X_{1}^{\bullet t}=\frac{k}{t} \biggr)= \biggl( \frac{t}{t+1}
\biggr)^k \frac{1}{t+1 },\qquad k=0,1,\dots.
\end{equation}
Thus, $t\sum_{i=1}^{M}X_{i}^{\bullet t}$ is a geometric sum of
geometric distributions with parameter $r=(1+t)^{-1}$. It is easy to
check (use LS transforms, for instance) that the distribution of such a
random variable is a mixture of the degenerate distribution at 0 (with
probability $p$) and a geometric distribution with parameter
$p^{*}=1-(1-r)(1-(1-p)r)^{-1}=1-t(t+p)^{-1}$, so that for each $k\in
\mathbb{N},$
%
\begin{eqnarray}\label{midisc}
L_{t}^{*}F \biggl(\frac{k}{t} \biggr) &=&
P\Biggl(\sum_{i=1}^{M}X_{i}^{\bullet t} \leq\frac
{k}{t}\Biggr) \nonumber\\[-8pt]\\[-8pt]
&=& p + (1 -
p)\bigl(1-(1-p^{*})^{k+1} \bigr) = 1 - (1 - p) \biggl( \frac{t}{t+p}
\biggr)^{k+1} .\nonumber
\end{eqnarray}
Note that the first equality in (\ref{midisc}) follows by recalling
(\ref{dfdfdi}) and noting that $(\sum_{i=1}^{M}X_{i} )^{\bullet t}$
has the same distribution as $\sum_{i=1}^{M}X_{i}^{\bullet t}$ (see
\cite{saerro}, Proposition 2.1). Actually, a more natural way (in this
case) to compute (\ref{midisc}) is to evaluate the LS transform of
$(\sum_{i=1}^{M}X_{i} )^{\bullet t}$ and then apply (\ref{discre}) and
(\ref{dfdfdi}). However, the previous computations enable easier
comparisons with the following method. In fact, one of the most obvious
(and widely used) methods to discretize the summands in (\ref{ransum})
is by a rounding method. For instance, a rounding down method (we round
$X_{i}$ to $[tX_{i}]t^{-1}$) yields
%
\begin{equation}\label{discal}
P \biggl(\frac{[tX_{1}]}{t}=\frac{k}{t} \biggr)=
P \biggl(\frac{k}{t}\leq X_{1}<\frac{k+1}{t} \biggr)=\mathrm{e}^{-k/t}(1-
\mathrm{e}^{-1/t}),\qquad k\in\mathbb{N}.
\end{equation}

In this case, $\sum_{i=1}^{M}[tX_{i}]$ is a geometric sum of geometric
distributions with parameter $r'=1-
\mathrm{e}^{-1/t}$. We denote by $R_{t}F$ the distribution function of $\sum
_{i=1}^{M}\frac{[tX_{i}]}{t}$. Using the same arguments as for (\ref
{midisc}), we obtain for each $k\in\mathbb{N}$ that
%
\begin{equation} \label{rodisc}
R_{t}F \biggl(\frac{k}{t} \biggr)=P\Biggl(\sum_{i=1}^{M}\frac
{[tX_{i}]}{t}\leq
\frac{k}{t}\Biggr)=1-(1-p) \biggl(\frac{\mathrm{e}^{-1/t}}{1-(1-p)(1-\mathrm{e}^{-1/t})}
\biggr)^{k+1} .
\end{equation}

Finally, it would be interesting to compare the previous
`discretization methods' with a `transform method.' To this end, we
consider the Laplace transform of $F$ in (\ref{testfi}) (instead of its
LS transform), that is,
\[
w_{F}(\theta)=\int_{0}^{\infty}\mathrm{e}^{-\theta u}F(u)\,\mathrm{d}u=\frac{1}{\theta
}-\frac{1-p}{\theta+p},\qquad \theta> 0,
\]
and apply the Post--Widder inversion formula (see \cite{feanin}, page
233), defined
for $t\in\mathbb{N}^{*}$ as
\[
W_{t}F(x)=\frac{(-1)^{t-1}}{(t-1)!} \biggl(\frac{t}{x} \biggr)^{t}w_{F}^{(t-1)}
\biggl(\frac{t}{x} \biggr)=1-\frac{(1-p)t^{t}}{(px+t)^{t}},\qquad x\geq0.
\]
In Table~\ref{table1} (computations with MATLAB) we consider a `rough'
discretization interval ($t=5$), a small $p$ ($p=0.1$) and present, for
different $x=k/5$, the exact values of $F$ (column 2), the $L_{t}^{*}$
approximation (column 3), the `rounding down' discretization (column 4)
and the Post--Widder inversion (column 5).

\begin{table}
\tablewidth=240pt
\caption{Comparison of different approximation methods for (\protect\ref{testfi})}\label{table1}
\begin{tabular*}{240pt}{@{\extracolsep{\fill}}lllll@{}}
\hline
$x=\frac{k}{5} $& $F (\frac{k}{5} )$ & $L_{5}^{*}F (\frac{k}{5} )$ &
$R_{5}F (\frac{k}{5} )$ & $ W_{5}F (\frac{k}{5} ) $\\
\hline
$0=\frac{0}{5}$& 0.1000 & 0.1176 & 0.1195 & 0.1000\\[2pt]
$1=\frac{5}{5}$& 0.1856 & 0.2008 & 0.2108 & 0.1848 \\[2pt]
$5=\frac{25}{5}$& 0.4541 & 0.4622 & 0.4907 & 0.4412 \\[2pt]
$10=\frac{50}{5}$& 0.6689 & 0.6722 & 0.7054 & 0.6383 \\[2pt]
$15=\frac{75}{5}$& 0.7992 & 0.8002 & 0.8296 & 0.7576\\[2pt]
$20=\frac{100}{5}$& 0.8782 & 0.8782 & 0.9014 & 0.8327\\[2pt]
$30=\frac{150}{5}$& 0.9552 & 0.9548 & 0.9670 & 0.9142\\[2pt]
$40=\frac{200}{5}$& 0.9835 & 0.9832 & 0.9890 & 0.9524\\
\hline
\end{tabular*}
\end{table}

As we can see in Table~\ref{table1}, $L_{5}^{*}F$ provides a better
approximation than $R_{5}F$. The intuitive explanation of this fact is
that, when approximating $\sum_{i=1}^{M}X_{i}$ by $\sum
_{i=1}^{M}X_{i}^{\bullet t}$, the error in the approximation can be controlled
`uniformly,' regardless of the distribution of $M$ (see \cite{saerro},
Theorem 4.3). This effect is obvious when we choose $M$ with a large
expected value (our choice of a small $p$ is for this reason -- for
larger values of $p$ checked, $L_{5}^{*}F$ is also better, but the
difference is less appreciable). However, if we compare the
approximations $L_{5}^{*}F$ and $W_{5}F$, we see that the last one is
better for small values, whereas the first one is better for large
values. To explain this fact, it is interesting to point out that
$W_{t}F$, like $L_{t}^{*}F$, admits the following well-known
representation. For a function $g$ defined on $[0,\infty)$, we can
write, as in (\ref{rlftlt}) (see \cite{feanin}, pages 220, 223),
%
\begin{equation}\label{operwi}
W_{t}g(x)=Eg \biggl(x\frac{S(t)}{t} \biggr),\qquad x>0.
\end{equation}
Note that the mean of the `random points' defining $W_{t}$ in (\ref{operwi}) is $E(xt^{-1}S(t))=x$,
whereas for $L_{t}^{*}$ in (\ref{rlftlt}), we have
$E(t^{-1}S([tx]+1))=t^{-1}([tx]+1)$. This means that $W_{t}$ is
centered at~$x$, whereas $L_{t}^{*}$ is `biased'. The benefits of this
property for $W_{t}$ are observed at small values in Table~\ref{table1}. However,
we have $\operatorname{Var}(xt^{-1}S(t))=t^{-1}x^{2}$, whereas
$\operatorname{Var}(t^{-1}S([tx]+1))=t^{-2}([tx]+1)$, the latter being of order
$t^{-1}x$, as $t\rightarrow\infty$. The greater variability of the
random variables defining $W_{t}$ for greater values of $x$ produces an
undesired effect in the approximation.

\begin{Table}
\caption{Comparison of $M_{5}^{[2]}$ in (\protect\ref{dftile}) with $
G_{5}^{[2]}$ in (\protect\ref{dfacpw})} \label{table2}
\begin{tabular}{@{}llllll@{}}
\hline
$x=\frac{k}{5} $& $F (\frac{k}{5} )$& $L_{5}^{*}F (\frac{k-1}{5} )$ &
$L_{10}^{*}F (\frac{2k-1}{10} )$ & $M_{5}^{[2]}F (\frac{k}{5} )$&$
G_{5}^{[2]}F (\frac{k}{5} ) $\\
\hline
$1=\frac{5}{5}$& 0.1856 & 0.1848 & 0.1852 & 0.1856 & 0.1856 \\[2pt]
$5=\frac{25}{5}$& 0.4541 & 0.4514 & 0.4528 & 0.4541 & 0.4538 \\[2pt]
$10=\frac{50}{5}$& 0.6689 & 0.6656 & 0.6673 & 0.6689 & 0.6677\\[2pt]
$15=\frac{75}{5}$& 0.7992 & 0.7962 & 0.7977 & 0.7992 & 0.7975\\[2pt]
$20=\frac{100}{5}$& 0.8782 & 0.8758 & 0.8770 & 0.8782 & 0.8766\\[2pt]
$30=\frac{150}{5}$& 0.9552 & 0.9538 & 0.9545 & 0.9552 & 0.9553\\[2pt]
$40=\frac{200}{5}$& 0.9835 & 0.9829 & 0.9832 & 0.9835 & 0.9854\\
\hline
\end{tabular}
\end{Table}

We now show the improvement in the approximation which occurs when
using $M_{t}^{[2]}$, as defined in (\ref{dftile}), instead of
$L_{t}^{*}$. In Table \ref{table2} below ($t=5$), we compare
$M_{t}^{[2]}F$ (column 5) with Richardson extrapolation for $W_{t}F$
(or Stehfest enhancement of order two for the Post--Widder formula --
see \cite{abwhla}, page 40), that is,
%
\begin{equation}\label{dfacpw}
G_{t}^{[2]}F(x):=2W_{2t}F(x)-W_{t}F(x),\qquad x>0.
\end{equation}
As we can see, $M_{5}^{[2]}F$ provides us with an exact value up to a
four decimal places, whereas $G_{5}^{[2]}F$ does not achieve this accuracy.
\end{example}

\section{Error bounds for the approximation}\label{seerro}

Let $g\in\mathcal{D}$, as defined in (\ref{dfclas}). Our first aim is
to give bounds of $\|L_{t}^{[2]}g-g\|$ in terms of
$\|x^{2}g^{iv}(x)\|$. To this end, we will use the following as `test function':
%
\begin{equation}\label{testfu}
\phi(x)= \cases{
0, & \quad if $x=0$, \vspace*{2pt}\cr
\displaystyle\frac{x^{2}}{2} \biggl(\frac{3}{2}- \log(x) \biggr), & \quad
otherwise. }
\end{equation}

Observe that $\phi\in\mathcal{D}$. In fact, by elementary calculus,
%
\begin{equation}\label{deriph}
\phi'(x)=x(1-\log x), \qquad
\phi''(x)=-\log x,\qquad
\phi'''(x)=-\frac{1}{x} \quad\mbox{and}\quad
\phi^{iv}(x)=\frac{1}{x^{2}}.
\end{equation}
In the
next lemma, we make an explicit computation of $L_{t}\phi(x)$ in
terms of the $\Psi$ (or digamma) function. This function
is defined as (see \cite{abstha}, page 258)
%
\begin{equation}\label{dfpsif}
\Psi(x):=\frac{\mathrm{d}}{\mathrm{d}x}
\log(\Gamma(x))=\frac{1}{\Gamma(x)}\int_{0}^{\infty}\log u \mathrm{e}^{-u}
u^{x-1}\,\mathrm{d}u, \qquad x>0,
\end{equation}
and, therefore, using
the last equality, we have the following probabilistic expression of
the psi function in terms of the gamma process:
%
\begin{equation}\label{psiope}
\Psi(x)=E\log S(x),\qquad x> 0.
\end{equation}

We will use the following property of this
function (see \cite{abstha}, page 258):
\begin{equation}\label{indpsi}
\Psi(x+1)=\frac{1}{x}+\Psi(x).
\end{equation}

\begin{lemma}\label{lephps}
Let $\phi$ be defined as in \textup{(\ref{testfu})}
and let
$L_{t}, t>0$, be defined as in \textup{(\ref{operbue})}. We have
that
%
\begin{equation}\label{operph}
L_{t}\phi(x)=\frac{1}{2t^{2}} \biggl(\frac{3(tx)^{2}}{2}-\frac{tx}{2}-1+tx(tx+1)\bigl(-\Psi
(tx)+\log(t)\bigr) \biggr),\qquad x>0.
\end{equation}
\end{lemma}

\begin{pf}
Let $t>0$ and $x>0$ be fixed. First, using elementary
calculus, (\ref{operbue}) and (\ref{indpsi}), we can write
\begin{eqnarray}
L_{t}\phi(x)&=&E
\frac{S(tx)^{2}}{2t^{2}} \biggl(\frac{3}{2}-
\log\biggl(\frac{S(tx)}{t} \biggr) \biggr)\nonumber\\
&=&\frac{1}{2t^2}\frac{1}{\Gamma(tx)}
\int_{0}^{\infty}u^{2}
\biggl(\frac{3}{2}-\log\biggl(\frac{u}{t} \biggr)
\biggr)\mathrm{e}^{-u}u^{tx-1}\,\mathrm{d}u\nonumber\\[-8pt]\\[-8pt]
&=&
\frac{(tx)(tx+1)}{2t^2}\frac{1}{\Gamma(tx+2)}\int_{0}^{\infty}
\biggl(\frac{3}{2}-\log\biggl(\frac{u}{t} \biggr) \biggr)\mathrm{e}^{-u}u^{tx+1}\,\mathrm{d}u\nonumber\\
&=&
\frac{(tx)(tx+1)}{2t^{2}} \biggl(\frac{3}{2}-
E\log\biggl(\frac{S(tx+2)}{t} \biggr) \biggr).\nonumber
\end{eqnarray}
Therefore,
using (\ref{psiope}), we can write
%
\begin{equation}\label{letes1}
L_{t}\phi(x)=\frac{(tx)(tx+1)}{2t^{2}} \biggl( \frac{3}{2}-
\Psi(tx+2)+\log(t) \biggr).
\end{equation}
Now, using
(\ref{indpsi}) twice, we have
%
\begin{equation}\Psi(tx+2)=\frac{2(tx)+1}{tx(tx+1)}+\Psi(tx).\label{letes2}
\end{equation}
By (\ref{letes1}), (\ref{letes2}), we obtain
\[
L_{t}\phi(x)=\frac{(tx)(tx+1)}{2t^{2}} \biggl(\frac{3}{2}-\frac
{2(tx)+1}{tx(tx+1)}-
\Psi(tx)+\log(t) \biggr).
\]
The result follows using elementary
algebra in the expression above.
\end{pf}

In the next lemma, we will study the approximation properties of
$L_{t}\phi$ to $\phi$. We will make use of the following
inequalities for the psi function:
%
\begin{eqnarray}\label{firpsi}
&\displaystyle \frac{1}{2x}\leq\log(x) - \Psi(x) \leq\frac{1
}{x},\qquad x>0; &
\\
\label{secpsi}
&\displaystyle \log(x)- \Psi(x) -\frac{1 }{2x} \leq\frac{1}{
12x^{2}},\qquad x>0.&
\end{eqnarray}
We can find (\ref{firpsi})
in \cite{alonso}, page 374, whereas (\ref{secpsi}) is an
immediate consequence of the fact that the function
\[
\Psi(x)-\log(x) +\frac{1 }{2x}+\frac{1}{
12x^{2}}
\]
is completely monotonic (see \cite{qicuso}, page 304) and
thus non-negative.

\begin{lemma}\label{letest}
Let $\phi$ be as defined in \textup{(\ref{testfu})} and
let $L_{t},
t>0$, be as defined in
\textup{(\ref{operbue})}. We have
%
\begin{equation}\label{letesb}
\biggl\|L_{t}\phi(x)-\phi(x)+\frac{x \log
x}{2t}+\frac{1}{3t^{2}}\biggr\|\leq
\frac{3}{8t^{2}}.
\end{equation}
\end{lemma}

\begin{pf}
Let $x>0$ and $t>0$ be fixed. First of all, we can write
%
\begin{equation}\label{funcph}
\phi(x)=\frac{1}{2t^{2}} \biggl(\frac
{3(tx)^{2}}{2}-(tx)^{2}\log(tx)+(tx)^{2}\log
(t) \biggr).
\end{equation}
On the other hand,
%
\begin{equation}\label{derip1}
\frac{x \log
x}{2t}+\frac{1}{3t^{2}}= \frac{1}{2t^{2}} \biggl((tx )\log tx-(tx)\log
t+\frac{2}{3} \biggr).
\end{equation}
Therefore, using
Lemma \ref{lephps}, (\ref{funcph}) and (\ref{derip1}), we can write
%
\begin{eqnarray}\label{normph1}
&&L_{t}\phi(x)-\phi(x)+\frac{x \log
x}{2t}+\frac{1}{3t^{2}}\nonumber\\
&&\quad=\frac{1}{2t^{2}} \biggl(-\frac
{tx}{2}-1-(tx)^{2}\Psi(tx)-
(tx)\Psi(tx)+(tx)^{2}\log(tx)+(tx)\log(tx)+\frac{2}{3} \biggr)\qquad
\\
&&\quad=\frac{1}{2t^{2}} \biggl((tx)^{2} \biggl(\log(tx)-\Psi(tx)-\frac
{1}{2(tx)} \biggr)+tx\bigl(\log(tx)-\Psi(tx)\bigr)-\frac{1}{3} \biggr).\nonumber
\end{eqnarray}
By (\ref{firpsi}), we have that $1/2\leq x(\log(x)-\Psi(x))\leq1,
x>0$, and thus
%
\begin{equation}\label{normph2}
\tfrac{1}{6}\leq
tx\bigl(\log(tx)-\Psi(tx)\bigr)-\tfrac{1}{3}\leq
\tfrac{2}{3}.
\end{equation}
Thus, using
(\ref{normph1}), (\ref{normph2}) and (\ref{secpsi}), we obtain
(\ref{letesb}).
\end{pf}

We are now in a position to state the following.

\begin{theorem}\label{teacce}
Let $g \in\mathcal{D}$, as defined in \textup{(\ref
{dfclas})} and
let $L^{[2]}_{t}$, $t>0$, be
as defined in \textup{(\ref{acgamo})}. We have
\[
\bigl|L^{[2]}_{t}g(x)-g(x)\bigr|\leq
\frac{1}{6t^{2}}\|xg'''(x)\|+\frac{9}{16t^{2}}\|x^{2}g^{iv}(x)\|
.
\]
\end{theorem}

\begin{pf}
We will first see that $g \in\mathcal{D}$ implies
that
%
\begin{equation} \|xg'''(x)\|\leq
\|x^{2}g^{iv}(x)\|<\infty.\label{bouthr}
\end{equation}
To begin with, the fact that
$\|x^{2}g^{iv}(x)\|<\infty$ implies that $\lim
_{x\rightarrow\infty}x^{1+\alpha}g^{iv}(x)=0$ for all $0<\alpha<1$. By
L'H\^{o}pital's rule, we also have that
$\lim_{x\rightarrow\infty}x^{\alpha}g'''(x)=0$, thus
concluding that
$\lim_{x\rightarrow\infty}g'''(x)=0$. Using this fact, we can write
\[
g'''(x)=-\int_{x}^{\infty}g^{iv}(u)\,\mathrm{d}u ,
\]
which implies easily (\ref{bouthr}) as
\[
|xg'''(x)|\leq x\int_{x}^{\infty}\frac
{|u^{2}g^{iv}(u)|}{u^{2}}\,\mathrm{d}u\leq
\|x^{2}g^{iv}(x)\| .
\]

Now, let $t>0$ and let $L_{t}$ be as
in (\ref{operbue}). As a previous step, we will prove that
%
\begin{equation}\label{teprim}
\biggl|L_{t}g(x)-g(x)-\frac{xg''(x)}{2t}-\frac
{xg'''(x)}{3t^{2}}\biggr|\leq
\frac{3}{8t^{2}}\|x^{2}g^{iv}(x)\|,\qquad
x>0.
\end{equation}
To this end, let $x>0$. Using a
Taylor series expansion of the random point $u=S(tx)/t$ around
$x$ and taking into account that $E(S(x)-x)=0$, $E(S(x)-x)^2=x$ and
$E(S(x)-x)^3=2x$, we can write
\begin{eqnarray}\label{taylor}
L_{t}g(x)-g(x)&=&Eg \biggl(\frac{S(tx)}{t} \biggr)-g(x)\nonumber\\
&=&\frac{E(S(tx)-tx)^2}{2t^{2}}g''(x)+\frac{E(S(tx)-tx)^3}{6t^{3}}g'''(x)\nonumber\\[-8pt]\\[-8pt]
&&{}+\frac{1}{6}E\int_{x}^{S(tx)/t} g^{iv}(\theta
)\biggl(\frac
{S(tx)}{t}-\theta\biggr)^{3}\,\mathrm{d}\theta\nonumber\\
&=&\frac{xg''(x)}{2t}+\frac{xg'''(x)}{3t^{2}}+\frac{1}{6}E\int
_{x}^{S(tx)/t}
g^{iv}(\theta)\biggl(\frac{S(tx)}{t}-\theta\biggr)^{3}\,\mathrm{d}\theta.\nonumber
\end{eqnarray}
Then, using (\ref{taylor}), we get the bound
\begin{eqnarray}\label{taybon}
&&\biggl|L_{t}g(x)-g(x)-\frac{xg''(x)}{2t}-\frac{xg'''(x)}{3t^{2}} \biggr|\nonumber\\
&&\quad=\frac
{1}{6} \biggl|E\int_{x}^{S(tx)/t} g^{iv}(\theta) \biggl(\frac
{S(tx)}{t}-\theta\biggr)^{3}\,\mathrm{d}\theta\biggr|
\nonumber\\[-8pt]\\[-8pt]
&&\quad\leq
\frac{\|x^{2}g^{iv}(x)\|}{6}E\int_{\min(x,S(tx)/t)}^{\max
(x,S(tx)/t)} \biggl|\frac{S(tx)}{t}-\theta\biggr|^{3}\frac{1}{\theta
^{2}}\,\mathrm{d}\theta\nonumber
\\
&&\quad=
\frac{\|x^{2}g^{iv}(x)\|}{6}E\int_{x}^{S(tx)/t} \biggl(\frac
{S(tx)}{t}-\theta\biggr)^{3}\frac{1}{\theta^{2}}\,\mathrm{d}\theta.\nonumber
\end{eqnarray}
Let $\phi(\cdot)$ be as in
(\ref{testfu}). Using (\ref{taylor}) and (\ref{deriph}), we have
%
\begin{equation}\label{interm}
L_{t}\phi(x)-\phi(x)+\frac{x \log
x}{2t}+\frac{1}{3t^{2}}
=\frac{1}{6}E\int_{x}^{S(tx)/t} \biggl(\frac{S(tx)}{t}-\theta
\biggr)^{3}\frac{1}{\theta^{2}}\,\mathrm{d}\theta.
\end{equation}
Then, by (\ref{taybon}) and
(\ref{interm}), we can write
\[
\biggl|L_{t}g(x)-g(x)-\frac{xg''(x)}{2t}-\frac{xg'''(x)}{3t^{2}}\biggr|\leq
\|x^{2}g^{iv}(x)\|\biggl\|L_{t}\phi(x)-\phi(x)+\frac{x \log
x}{2t}+\frac{1}{3t^{2}}\biggr\| .
\]
Thus, (\ref{teprim}) follows by applying
Lemma \ref{letest}.

Observe that in (\ref{teprim}), the only term of order $1/t$ is the
one involving the second derivative. By means
of the operator $L^{[2]}_{t}$, as defined in (\ref{acgamo}), this term
is eliminated. In fact, using (\ref{teprim}), we have
\begin{eqnarray}
L^{[2]}_{t}g(x)-g(x)&=&2 \bigl(L_{2t}g(x)-g(x) \bigr)-\bigl(L_{t}g(x)-g(x)\bigr)
\nonumber\\
&=&2 \biggl(L_{2t}g(x)-g(x)-\frac{x}{4t}g''(x)-\frac{x}{12t^{2}}g'''(x)
\biggr)\nonumber\\[-8pt]\\[-8pt]
&&{}- \biggl(L_{t}g(x)-g(x)-\frac{x}{2t}g''(x)-\frac{x}{3t^{2}}g'''(x) \biggr)
-\frac{x}{6t^{2}}g'''(x)\nonumber\\
&\leq&
\frac{1}{6t^{2}}\|xg'''(x)\|+\frac{9}{16t^{2}}\|x^{2}g^{iv}(x)\|.\nonumber
\end{eqnarray}
This completes the proof of Theorem \ref{teacce}.
\end{pf}

Finally, in the next result, we study the approximation
properties of $M^{[2]}_{t}$.

\begin{theorem}\label{tetota}
Let $g\in\mathcal{D}_{1}$, as
defined in \textup{(\ref{dfcla2})} and let $M^{[2]}_{t}, t>0$, be as defined
in \textup{(\ref{a2gamo})}. We have
\[
\bigl\|M^{[2]}_{t}g-g\bigr\|\leq\frac{1}{8t^{2}}\|g''(x)\|
+\frac{1}{6t^{2}}\|xg'''(x)\|+\frac{9}{16t^{2}}\|x^{2}g^{iv}(x)\|.
\]
\end{theorem}

\begin{pf}
First, note that $g\in
\mathcal{D}_{1}$ implies that $\|xg'''(x)\|<\infty$, thanks to
(\ref{bouthr}). Now, let $t>0$ and $x>0$ be fixed. We write
\begin{eqnarray}\label{prteto}
M^{[2]}_{t}g(x)-g(x)&=&(tx-[tx])
\biggl(L^{[2]}_{t}g \biggl(\frac{[tx]+1}{t} \biggr)-g \biggl(\frac{[tx]+1}{t} \biggr) \biggr)\nonumber\\
&&{} +([tx]+1-tx)
\biggl(L^{[2]}_{t}g \biggl(\frac{[tx]}{t} \biggr)-g \biggl(\frac{[tx]}{t} \biggr) \biggr)\nonumber\\[-8pt]\\[-8pt]
&&{} +(tx-[tx]) \biggl(g \biggl(\frac{[tx]+1}{t} \biggr)-g(x) \biggr)\nonumber\\
&&{}+([tx]+1-tx) \biggl(g
\biggl(\frac{[tx]}{t} \biggr)-g(x) \biggr).\nonumber
\end{eqnarray}
Using the usual expansion
%
\begin{equation}\label{taysec}
|g(y)-g(x)-(y-x)g'(x)|\leq \frac{(y-x)^{2}}{2}\|g''\|
\end{equation}
and taking into account that
%
\begin{eqnarray}
&&(tx-[tx]) \biggl(g \biggl(\frac{[tx]+1}{t} \biggr)-g(x) \biggr)+([tx]+1-tx) \biggl(g
\biggl(\frac{[tx]}{t} \biggr)-g(x) \biggr)\nonumber\\
&&\quad=(tx-[tx]) \biggl(g \biggl(\frac{[tx]+1}{t} \biggr)-g(x)-\frac{[tx]+1-tx}{t}g'(x)
\biggr)\\
&&\qquad{}+([tx]+1-tx) \biggl(g \biggl(\frac{[tx]}{t} \biggr)
-g(x)-\frac{[tx]-tx}{t}g'(x) \biggr),\nonumber
\end{eqnarray}
we obtain from the
above expression and (\ref{taysec})
that
%
\begin{eqnarray}\label{bouthi}
&&\biggl|(tx-[tx]) \biggl(g \biggl(\frac{[tx]+1}{t} \biggr)-g(x) \biggr)+([tx]+1-tx) \biggl(g
\biggl(\frac{[tx]}{t} \biggr)-g(x) \biggr) \biggr|
\nonumber\\
&&\quad\leq\biggl((tx-[tx])\frac
{([tx]+1-tx)^{2}}{2t^{2}}+([tx]+1-tx)\frac{([tx]-tx)^{2}}{2t^{2}} \biggr)\|
g''\|\\
&&\quad=\frac{(tx-[tx])([tx]+1-tx)}{2t^{2}}\|g''\|\leq
\frac{1}{8t^{2}}\|g''\|,\nonumber
\end{eqnarray}
the last
inequality holding since for each $k\in\mathbb{N}$, the supremum of
$(u-k)(k+1-u),
k\leq u\leq k+1$, is attained at $u=k+1/2$. On the other hand,
taking into account Theorem \ref{teacce}, we have
%
\begin{eqnarray}\label{boufis}
&&\biggl|(tx-[tx]) \biggl(L^{[2]}_{t}g \biggl(\frac{[tx]+1}{t} \biggr)-g \biggl(\frac
{[tx]+1}{t} \biggr) \biggr)\nonumber\\
&&\qquad{}+([tx]+1-tx)
\biggl(L^{[2]}_{t}g \biggl(\frac{[tx]}{t} \biggr)-g \biggl(\frac{[tx]}{t} \biggr) \biggr)
\biggr|\\
&&\quad\leq \bigl\|L^{[2]}_{t}g-g\bigr\| \leq
\frac{1}{6t^{2}}\|xg'''(x)\|+\frac{9}{16t^{2}}\|x^{2}g^{iv}(x)\|.\nonumber
\end{eqnarray}
The result follows by (\ref{prteto}), (\ref{bouthi}) and
(\ref{boufis}).
\end{pf}

\section{Application to gamma distributions}\label{seapp1}

In this section, we will study the case of gamma distributions, that
is, distributions with density functions as given in (\ref{gamden}). It
is not hard
to see that these distributions are in the class $\mathcal{D}_{1}$, for
a shape parameter $p=1$ or $p\geq2$, and, therefore, we are a
position of apply Theorem \ref{tetota}. The aim of this section is
to show that, in fact, the bounds in this theorem can be uniformly
bounded on the shape parameter, which will be an advantage when
dealing with mixtures of these distributions. From now on, we
define
%
\begin{equation}\label{gamde1}
f_{p}(x):= \cases{
\displaystyle\frac{\mathrm{e}^{-x}x^{p-1}}{\Gamma(p)},& \quad$x>0$, if $ p\in \mathbb{R}
\setminus\{0,-1,-2,\dots\}$,\vspace*{2pt}\cr
0, &\quad $x>0$, if $p\in\{0,-1,-2,\dots\}$. }
\end{equation}
The `odd' definition of $f_{p}$ for $ p\in\{0,-1,-2,\dots\}$ is for
notational convenience in (\ref{ledeg1}). For \mbox{$p>0$}, the
function above is the density of a gamma random variable as in
(\ref{gamden}), with scale parameter $a=1$. Results for another
scale parameter will follow by a change of scale (see Proposition~\ref{presca} below). First, we
will consider the case $p=1$, that is, an exponential random
variable. As the distribution function of this random variable presents
no computational problems, it makes no sense to approximate it.
However, when we consider the problem of approximating a general
mixture of Gamma distributions, the exponential distribution could
be a component.

\begin{lemma}\label{leboex}
Let $F(x)=1-\mathrm{e}^{-x}, x\geq0$. For $t>0$, let
$M_{t}^{[2]}F$ be as defined in \textup{(\ref{a2gamo})}. We have that
\[
\bigl\|M_{t}^{[2]}F-F\bigr\|\leq\biggl(\frac{1}{8}+\frac{1}{6 \mathrm{e}}+\frac{9}{4\mathrm{e}^{2}}
\biggr)\frac{1}{t^{2}}.
\]
\end{lemma}

\begin{pf}
First of all, note that $|F^{(k)}(x)|=\mathrm{e}^{-x}$ and that
$\sup_{x\geq0}x^{k}\mathrm{e}^{-x}=k^{k}\mathrm{e}^{-k},
k=1,2,\ldots.$ Thus, we have
%
\begin{equation}\label{bouexp}
\|F''\|=1, \qquad \|xF'''(x)\|=\mathrm{e}^{-1}
\quad\mbox{and}\quad
\|x^{2}F^{iv}(x)\|=2^{2}\mathrm{e}^{-2}.
\end{equation}
The
conclusion follows by taking into account Theorem \ref{tetota}.
\end{pf}

We will now deal with the case $p\geq2$ in (\ref{gamde1}). The
two following lemmas will be useful in order to bound the derivatives of
this density. For the sake of brevity, they are stated without proof
(only elementary calculus is required). For the proofs, we refer the
interested reader to \cite{saunif}, a preliminary version of this paper
(available online).

\begin{lemma}\label{ledega}
Let $f_{p}(\cdot)$, $p >0$, be as
defined in \textup{(\ref{gamde1})}. We have, for all $n\in\mathbb{N}$,
\begin{eqnarray}\label{ledeg1}
\frac{\mathrm{d}^n}{\mathrm{d}x^n}f_{p}(x)&=&\frac{\mathrm{e}^{-x}x^{p-n-1}}{\Gamma(p)}\sum_{i=0}^{n}\pmatrix{n \cr i}(-1)^{i} \Biggl(\prod_{j=1}^{n-i}(p-j)
\Biggr)x^{i}\nonumber\\[-8pt]\\[-8pt]
&=& \sum_{i=0}^{n}\pmatrix{n \cr i}(-1)^{i}
f_{p-n+i}(x),\qquad x>0,\nonumber
\end{eqnarray}
in which
$\prod_{j=1}^{0}(p-j)=1$.
\end{lemma}

Next, we formulate a technical lemma in which we define
certain decreasing functions which will be used to bound the
weighted derivatives of $f_{p}$.

\begin{lemma}\label{ledecr}
We have:
\begin{longlist}[(iii)]
\item[(i)] the function
%
\begin{equation}\label{dfgep1}
g_{1}(p):=\frac{1}{\Gamma(p)}\mathrm{e}^{-(p-1)}(p-1)^{p-1},\qquad
p> 1
\end{equation}
$(g_{1}(1)=1 )$, is
decreasing in $p$;
\item[(ii)] the function
%
\begin{equation}\label{dfgep2}
g_{2}(p):=\frac{1}{\Gamma(p)}\mathrm{e}^{- (p-
1/2+1/2\sqrt{4p-3} )} \biggl(p-
\frac{1}{2}+\frac{1}{2}\sqrt{4p-3} \biggr)^{p-1/2},\qquad p \geq
1,
\end{equation}
is decreasing in $p$;
\item[(iii)] the function
%
\begin{equation}\label{dfgep3}
g_{3}(p):=\frac{1}{\Gamma(p)}\mathrm{e}^{-(p-1-\sqrt
{p-1})}\bigl(\sqrt
{p-1}-1\bigr)^{p-2}\bigl(\sqrt{p-1}\bigr)^{p-1},\qquad p> 2
\end{equation}
$ (g_{3}(2)=1 )$, is decreasing in $p$;
\item[(iv)] the function
%
\begin{equation}\label{dfgep4}
g_{4}(p):=\frac{1}{\Gamma(p)}\mathrm{e}^{-(p-\sqrt
{3p-2})}\bigl(p-\sqrt{3p-2}\bigr)^{p-2}\bigl(\sqrt{3p-2}-1\bigr)^{3},\qquad
p>2
\end{equation}
$ (g_{4}(2)=1 )$, is
decreasing in $p$.
\end{longlist}
\end{lemma}

In the following result, we get bounds of the quantities
required in Theorem \ref{tetota}, depending on the shape parameter
$p$, but also decreasing on $p$.

\begin{lemma}\label{leboga}
Let $f_{p}$ be as in \textup{(\ref{gamde1})} and $g_{i},
i=1,2,3,4$, be as in Lemma \textup{\ref{ledecr}}. We have:
\begin{eqnarray*}
\mathrm{(i)}&\hspace*{6pt}&
\sup_{x\geq
0}|f_{p}(x)|=g_{1}(p),\qquad p\geq1;\\
\mathrm{(ii)}&\hspace*{6pt}&
\sup_{x\geq
0}|xf'_{p}(x)|=g_{2}(p),\qquad p\geq1;\\
\mathrm{(iii)}&\hspace*{6pt}&
\sup_{x\geq
0}|f_{p}'(x)|=g_{3}(p),\qquad p\geq2;\\
\mathrm{(iv)}&\hspace*{6pt}&
\sup_{x\geq 0}|xf_{p}''(x)|\leq\max\{g_{1}(p-1),g_{2}(p-1)\},\qquad
p\geq2;\\
\mathrm{(v)}&\hspace*{6pt}&
\sup_{x\geq
0}|x^{2}f_{p}'''(x)|\leq g_{4}(p)+3g_{2}(p-1)+g_{1}(p-1),\qquad p\geq
2.
\end{eqnarray*}
\end{lemma}

\begin{pf}
To show part (i), it is clear that, for $p\geq1$,
\[
\sup_{x\geq
0}f_{p}(x)=f_{p}(p-1)=\frac{\mathrm{e}^{-(p-1)}(p-1)^{p-1}}{\Gamma(p)}
\]
and
(i) follows by recalling (\ref{dfgep1}). To show part (ii), we have (see
\cite{saerro}, Remark 3.2 and Lemma~5.2)
%
\begin{equation}\label{partb2}
\sup_{x\geq0}|xf_{p}'(x)|=\frac{1}{\Gamma(p)} \biggl(p-
\frac{1}{2}+\frac{1}{2}\sqrt{4p-3} \biggr)^{p-1/2}\mathrm{e}^{-p-
1/2+1/2\sqrt{4p-3}},\qquad p>1,
\end{equation}
and
(ii) follows by recalling (\ref{dfgep2}). To show part (iii), by
(\ref{ledeg1}), we have for $p\geq2$ that
%
\begin{eqnarray}\label{gamdem}
f_{p}'(x)&=&\frac{1}{\Gamma(p)}\mathrm{e}^{-x}x^{p-2}(p-1-x),\qquad
x>0,\\
\label{gamde2}
f_{p}''(x)&=&\frac{1}{\Gamma
(p)}\mathrm{e}^{-x}x^{p-3}\bigl((p-1)(p-2)-2(p-1)x+x^{2}\bigr),\qquad
x>0,
\end{eqnarray}
and it can be easily checked that
the zeros of $f''_{p}(x)$ are $p_{1}:=p-1-\sqrt{p-1}$ and
$p_{2}:=p-1+\sqrt{p-1}$. Therefore, $|f_{p}'(x)|$ must attain its
maximum value at either $p_{1}$ or $p_{2}$. Actually, $p_{1}$
corresponds to the maximum. To show that, we will see that
%
\begin{equation}\label{quotie}
\frac{f_{p}'(p_{1})}{|f_{p}'(p_{2})|}=\mathrm{e}^{2\sqrt{p-1}}
\biggl(\frac{\sqrt{p-1}-1}{\sqrt{p-1}+1} \biggr)^{p-2}\geq1,\qquad
p\geq2.
\end{equation}
To show the last inequality in
(\ref{quotie}), taking logarithms, we will prove that
%
\begin{equation}\label{quotie1}
r_{1}(p):= 2\sqrt{p-1} + (p-2) \bigl(\log\bigl(\sqrt{p-1}-1\bigr)
- \log
\bigl(\sqrt{p-1}+1\bigr) \bigr)\geq0,\qquad
p> 2.
\end{equation}
Define
\[
\rho_{1}(b):=\frac{2b}{b^2-1}+ \bigl(\log(b-1)-\log(b+1) \bigr),\qquad b>1.
\]
Note that
%
\begin{equation}\label{reros1}
r_{1}(p)=(p-2)\rho_{1}\bigl(\sqrt{p-1}\bigr),\qquad p>2.
\end{equation}
We will first prove that
%
\begin{equation}\label{quotie2}
\rho_{1}(b)\geq0,\qquad b>1.
\end{equation}
To show (\ref{quotie2}), it is
readily seen that $\rho_{1}'(b)=-4(b^2-1)^{-2}, b>1$, so that
$\rho_{1}$ is decreasing. As $\lim_{b\rightarrow
\infty}\rho_{1}(b)=0$, we have (\ref{quotie2}). This implies also
(\ref{quotie1}), recalling (\ref{reros1}). Therefore, we conclude
that
%
\begin{equation}\label{desco1}
\sup_{x>0}|f_{p}'(x)|=f_{p}'(p_{1})=
\frac{1}{\Gamma(p)}\mathrm{e}^{-(p-1-\sqrt{p-1})}\bigl(\sqrt{p-1}-1\bigr)^{p-2}\bigl(\sqrt
{p-1}\bigr)^{p-1},
\end{equation}
which, together with (\ref{dfgep3}), shows (iii).

To show part (iv), note that by using (\ref{ledeg1}), we can write
$f_{p}'(x)=f_{p-1}(x)-f_{p}(x)$ and, therefore,
%
\begin{equation}\label{partb1}
xf_{p}''(x)=xf_{p-1}'(x)-xf_{p}'(x),\qquad x>0, p\geq2.
\end{equation}
On the other hand, we see in (\ref{gamde2}) that $f_{p-1}'(x)$ and
$f_{p}'(x)$ have the same sign for $0<x<p-2$ and $p-1<x<\infty$ and,
therefore, using part (ii) and Lemma \ref{ledecr}(i), we can write
%
\begin{equation}\label{bound1}
\sup_{x\notin[p-2,p-1]}|xf_{p}''(x)|\leq
\max\bigl(g_{2}(p-1),g_{2}(p)\bigr)=g_{2}(p-1) .
\end{equation}

On the other hand, we have, by (\ref{gamde2}),
%
\begin{equation}\label{partbr}
xf_{p}''(x)=\frac{1}{\Gamma
(p)}\mathrm{e}^{-x}x^{p-2}\bigl((p-1)(p-2)-2(p-1)x+x^{2}\bigr).
\end{equation}
Using the above expression and taking
into account that, for $p-2\leq x\leq p-1$,
%
\begin{equation}\label{partb3}
\mathrm{e}^{-x}x^{(p-2)}\leq \mathrm{e}^{-p-2}(p-2)^{p-2}
\quad\mbox{and}\quad
|(p-1)(p-2)-2(p-1)x+x^{2}|=p-1,
\end{equation}
the last inequality holds as
$|(p-1)(p-2)-2(p-1)x+x^{2}|, p-2\leq x\leq p-1$, attains its maximum
value at $p-1$. From (\ref{partbr}) and (\ref{partb3}), we conclude
that
%
\begin{equation}\label{bound2}
\sup_{x \in[p-2,p-1]}|xf_{p}''(x)|\leq
\frac{1}{\Gamma(p)}\mathrm{e}^{-(p-2)}(p-2)^{p-2}(p-1)=g_{1}(p-1),
\end{equation}
where the last inequality follows by recalling (\ref{dfgep1}). Thus,
(\ref{bound1}) and (\ref{bound2})
conclude the proof of part~(iv). To show part (v), let $p\geq2$. First,
we have, by (\ref{ledeg1}),
%
\begin{eqnarray}\label{gamde3}
f_{p}'''(x)&=&f_{p-3}(x)-3f_{p-2}(x)+3f_{p-1}(x)-f_{p}(x)\nonumber
\\
&=&\frac{\mathrm{e}^{-x}x^{p-4}}{\Gamma
(p)}\bigl((p-1)(p-2)(p-3)-3(p-1)(p-2)x+3(p-1)x^{2}-x^{3}\bigr)\qquad\\
&=&\frac{\mathrm{e}^{-x}x^{p-4}}{\Gamma(p)}\bigl((p-1-x)^{3} +
3(p-1)\bigl(x-(p-2)\bigr) - (p-1)\bigr),\qquad
x>0.\nonumber
\end{eqnarray}

Therefore, if we call
\[
h_{p}(x):=\frac{\mathrm{e}^{-x}x^{p-2}}{\Gamma(p)}(p-1-x)^{3},\qquad x> 0,
\]
we have, recalling (\ref{gamdem}),
\begin{eqnarray}\label{botird}
x^{2}f_{p}'''(x)&=&\frac{\mathrm{e}^{-x}x^{p-2}}{\Gamma
(p)}\bigl((p-1-x)^{3}-3(p-1)\bigl(x-(p-2)\bigr)-(p-1)\bigr)\nonumber\\[-8pt]\\[-8pt]
&=&h_{p}(x)+3xf'_{p-1}(x)-f_{p-1}(x),\qquad x\geq\nonumber
0.
\end{eqnarray}
We will firstly see that
%
\begin{equation}\label{maxige}
\sup_{x\geq0}|h_{p}(x)|=g_{4}(p)
\end{equation}
with $g_{4}(\cdot)$ as defined in (\ref{dfgep4}). Note that
\[
h'_{p}(x)=\frac{\mathrm{e}^{-x}x^{p-3}}{\Gamma(p)}(p-1-x)^{2}\bigl(x^{2}-2p
x+(p-1)(p-2)\bigr),\qquad x>0 .
\]
The maximum value of $|h_{p}|$ will be
attained at the roots of the last polynomials, being
$p_{1}:=p+\sqrt{3p-2}$ and $p_{2}:=p-\sqrt{3p-2}$. To check which
value attains the maximum, define $u:=\sqrt{3p-2}$. Note that
$p_{1}=(u+1)(u+2)/3$ and $p_{2}=(u-1)(u-2)/3$. Then, with this
notation, we will prove that
%
\begin{equation}\label{maxroo}
\frac{|h_{p}(p_{2})|}{|h_{p}(p_{1})|}=\mathrm{e}^{2u} \biggl(\frac
{(u-1)(u-2)}{(u+1)(u+2)} \biggr)^{(u^2-4)/3}
\biggl(\frac{u-1}{u+1} \biggr)^{3}\geq1,\qquad
u>2.
\end{equation}
To show the last inequality in
(\ref{maxroo}), taking logarithms, we will show that
%
\begin{equation}\label{lnmaxr}
\rho_{2}(u):=2u+\frac{u^2-4}{3}\log\biggl(\frac
{(u-1)(u-2)}{(u+1)(u+2)} \biggr)+3\log\biggl(\frac{u-1}{u+1}
\biggr)\geq0,\qquad u>2.
\end{equation}
Note that
\begin{eqnarray*}
\rho_{2}'(u)&=&2+\frac{2u}{3}\log\biggl(\frac
{(u-1)(u-2)}{(u+1)(u+2)} \biggr)+\frac{u^{2}-4}{3}
\biggl(\frac{1}{u-1}+\frac{1}{u-2}-\frac{1}{u+1}-\frac{1}{u+2} \biggr)\\&&{}+3
\biggl(\frac
{1}{u-1}-\frac{1}{u+1} \biggr)\\
&=&\frac{4u^{2}}{u^2-1}+\frac{2u}{3}\log
\biggl(\frac
{(u-1)(u-2)}{(u+1)(u+2)} \biggr),\qquad
u>2 .
\end{eqnarray*}
We will show that $\rho_{2}'(u)\leq0$, $u>2$. In
fact,
\[
\frac{\mathrm{d}}{\mathrm{d}u}\frac{3}{2u}\rho_{2}'(u)=\frac
{36}{(u+1)^2(u-1)^2(u^2-4)^2}\geq0,\qquad
u>2 ,
\]
and then $3(2u)^{-1}\rho_{2}'(u)$ is increasing. As
$\lim_{u\rightarrow\infty}3(2u)^{-1}\rho_{2}'(u)=0$, we conclude
that $3(2u)^{-1}\times\rho_{2}'(u)\leq0$ and thus that $\rho_{2}'(u)\leq
0$. Therefore, $\rho_{2}(u)$ is decreasing. This, together with the
fact that $\lim_{u\rightarrow\infty}\rho_{2}(u)=0$, proves
(\ref{lnmaxr}) and therefore (\ref{maxroo}). Then,
$\|h_{p}\|=h_{p}(p_{2})=g_{4}(p)$, thus proving (\ref{maxige}). The
proof of part (iv) now follows easily by recalling (\ref{botird}) and
using (\ref{maxige}) and parts~(i) and (ii).
\end{pf}

As an immediate consequence of Theorem \ref{tetota} and Lemma
\ref{leboga}, we have the following corollary.

\begin{corollary}\label{coboga}
Let $F_{p}$ be a gamma distribution with shape
parameter $p\geq2$, that is, whose density
function is given by \textup{(\ref{gamde1})}. Let $M^{[2]}_{t}, t>0$,
be defined
as in \textup{(\ref{a2gamo})}. We have
\[
\bigl\|M^{[2]}_{t}F_{p}-F_{p}\bigr\|\leq
\biggl(\frac{17}{12}+\frac{27}{16\mathrm{e}} \biggr)\frac{1}{t^{2}}\approx
\frac{2.0375}{t^{2}}.
\]
\end{corollary}

\begin{pf}
Let $p \geq2$ be fixed. The result is an immediate
consequence of Theorem
\ref{tetota}, as $F_{p}'=f_{p}$, as defined in (\ref{gamde1}).
Therefore, by Lemma \ref{leboga}(iii) and Lemma \ref{ledecr}(ii), we
have that
%
\begin{equation}\label{cobog1}
\|F_{p}''\|=\|f_{p}'\|=g_{3}(p)\leq g_{3}(2)=1.
\end{equation}
On the other hand, we see that by Lemma \ref{ledecr}(i), we have that
%
\begin{equation}\label{cobogm}
g_{1}(p-1)\leq g_{1}(1)=1 \quad\mbox{and}\quad
g_{2}(p-1)\leq g_{2}(1)=\mathrm{e}^{-1},\qquad p\geq 2 .
\end{equation}
Thus, using the above inequalities and
Lemma \ref{leboga}(iv), we have
%
\begin{equation}\label{cobog2}
\|xF_{p}'''(x)\|=\|xf_{p}''(x)\|\leq1.
\end{equation}
Finally by Lemma \ref{leboga}(v), Lemma \ref{ledecr}(iv) and
(\ref{cobogm}), we have
%
\begin{equation}\label{cobog3}
\|x^{2}F_{p}^{iv}(x)\|=\|x^{2}f_{p}'''(x)\|\leq
g_{4}(2)+3g_{2}(1)+g_{1}(1)=2+3\mathrm{e}^{-1}.
\end{equation}
Using (\ref{cobog1}), (\ref{cobog2}), (\ref{cobog3}) and Theorem
\ref{tetota}, we obtain the result. This completes the proof of
Corollary \ref{coboga}.
\end{pf}

\section{Applications to mixtures of Erlang distributions and
phase-type distributions}\label{seapp2}

In this section we apply the results from the previous
section to mixtures of Erlang distributions and to random sums of
thereof. In order to undertake this study for an arbitrary scale parameter,
we need the following result which shows the behavior of $M_{t}^{[2]}F$ under
changes of scale.

\begin{proposition}\label{presca}
Let $X$ be a random variable with distribution function $F$. For a
given $c>0$, denote by $F^{c}$ the distribution function of $cX$. Let
$M_{t}^{[2]}F$ and $M_{t}^{[2]}F^{c}$, $t>0$, be the respective
approximations for $F$ and $F^{c}$, as defined in \textup{(\ref
{a2gamo})}. We
have that
%
\begin{equation}\label{presc1}
M_{t}^{[2]}F^{c}(x)=M_{ct}^{[2]}F(x/c),\qquad x\geq0.
\end{equation}
Therefore,
%
\begin{equation}\label{presc2}
\bigl\|M_{t}^{[2]}F^{c}-F^{c}\bigr\|=\bigl\|M_{ct}^{[2]}F-F\bigr\|.
\end{equation}
\end{proposition}

\begin{pf}
Let $t>0$ and $c>0$ be fixed. First, we will see that
%
\begin{equation}\label{prsc1l}
M_{t}^{[2]}F^{c} \biggl(\frac{k}{t} \biggr)=M_{ct}^{[2]}F \biggl(\frac
{k}{ct} \biggr),\qquad k\in\mathbb{N},
\end{equation}
and, therefore, (\ref{presc1}) is satisfied for points in the set
$k/t, k\in\mathbb{N}$. To this end, we use (\ref{dftil1}) and
(\ref{acgamo}), and take into account that
%
\begin{equation}\label{rescal}
F^{c}(x)=F(x/c),\qquad x\geq0,
\end{equation}
to
write, for all $k\in\mathbb{N}$,
\begin{eqnarray}
M_{t}^{[2]}F^{c} \biggl(\frac{k}{t} \biggr)&=&2EF^{c} \biggl(\frac{S(2k)}{2t} \biggr)-EF^{c} \biggl(\frac{S(k)}{t}
\biggr)\nonumber\\[-8pt]\\[-8pt]
&=&2EF
\biggl(\frac{S(2k)}{2ct} \biggr)-EF \biggl(\frac{S(k)}{ct} \biggr)=M_{ct}^{[2]}F \biggl(\frac
{k}{ct} \biggr),\nonumber
\end{eqnarray}
thus proving (\ref{prsc1l}). For a general $x>0$, we use
(\ref{a2gamo}) and (\ref{prsc1l}), to see that
\begin{eqnarray*}
M_{t}^{[2]}F^{c}(x)&=&(tx-[tx])M^{[2]}_{t}F^{c} \biggl(\frac
{[tx]+1}{t} \biggr)+([tx]+1-tx)M^{[2]}_{t}F^{c} \biggl(\frac{[tx]}{t} \biggr)
\\
&=&(tx-[tx])M^{[2]}_{ct}F \biggl(\frac{[tx]+1}{ct}
\biggr)+([tx]+1-tx)M^{[2]}_{ct}F \biggl(\frac{[tx]}{ct} \biggr)
=M_{ct}^{[2]}F \biggl(\frac{x}{c} \biggr),
\end{eqnarray*}
the last inequality being trivial as $tx=(ct)(x/c)$. This
concludes the proof of (\ref{presc1}). Finally, (\ref{presc2})
follows easily from (\ref{presc1}) and (\ref{rescal}), as we have
\[
\sup_{x>0}\bigl|M_{t}^{[2]}F^{c}(x)-F^{c}(x)\bigr|=\sup_{x>0}\bigl|
M_{ct}^{[2]}F(x/c)-F(x/c)\bigr|.
\]
This concludes the proof of Proposition
\ref{presca}.
\end{pf}

As an application of the results in the previous section, we will
consider the class of (possibly infinite) mixtures of Erlang
distributions recently studied by Willmot and Woo (see
\cite{wiwoon}). More specifically, let $F_{(a,j)} , a>0, j\in
\mathbb{N}^{*}$, be the distribution function corresponding to the density
$f_{(a,j)}$ given in (\ref{gamden}) (an Erlang $j$ distribution
with scale parameter $a$). We will consider a
finite number of scale parameters arranged in increasing order
($0<a_{1}<\cdots<a_{n}$) and a set of non-negative numbers $p_{ij},
i=1,\ldots, n, j=0,1,2,\ldots,$ such that $
\sum_{i=1}^{n}\sum_{j=1}^{\infty}p_{ij}=p\leq1$, and define the
class of
distribution functions $\mathcal{ME}(a_{1},\ldots, a_{n})$ given as
%
\begin{equation}\label{mixgav}
F(x)=(1-p)
+\sum_{i=1}^{n}\sum_{j=1}^{\infty}p_{ij}F_{a_{i},j}(x),\qquad x\geq0
\end{equation}
(we consider a slight modification of the class in \cite{wiwoon}, page
103, as we allow the point mass at 0 with probability $1-p$). Based on
\cite{wiwoon}, page 103, we can alternatively write (\ref{mixgav})
by using only the maximum of the scale parameters, that is,
%
\begin{equation}\label{mixgag}
F(x)=(1-p)+\sum_{j=1}^{\infty}p_{j}F_{a_{n},j}(x),\qquad
x\geq0.
\end{equation}
Moreover, the class
(\ref{mixgag}) is a wide class containing many of the distributions
considered in applied probability, such as (obviously) finite
mixtures of Erlang distributions, but also the class of phase-type distributions
(see Proposition \ref{cophas} below). Every random variable having a
representation as in (\ref{mixgav}) can be approximated by means of
$M_{t}^{[2]}$, as shown in the following result.

\begin{proposition}\label{procot}
Let $F$ be a distribution function of the form
$\mathcal{ME}(a_{1},\ldots,
a_{n})$, $0<a_{1}<\cdots<a_{n}$, as in \textup{(\ref{mixgav})}. Let
$M_{t}^{[2]}, t>0$, be defined as in \textup{(\ref{a2gamo})}. We have
%
\begin{equation}\label{cotmix}
\bigl\|M_{t}^{[2]}F-F\bigr\|\leq\biggl(\frac{17}{12}+\frac{27}{16\mathrm{e}} \biggr)\frac{ \sum
_{i=1}^{n}(\sum_{j=1}^{\infty}p_{ij})a_{i}^2}{t^2}.
\end{equation}
\end{proposition}

\begin{pf}
Let $t>0$ and $0<a_{1}<\cdots<a_{n}$ be fixed. The linearity
of $M_{t}^{[2]}$ yields
%
\begin{equation}\label{linear}
M_{t}^{[2]}F(x)=(1-p)+\sum_{i=1}^{n}\sum
_{j=1}^{\infty
}p_{ij}M_{t}^{[2]}F_{a_{i},j}(x),\qquad
x\geq0.
\end{equation}
By Corollary \ref{coboga}, we
can write, for a scale parameter 1,
%
\begin{equation}\label{lemtmj}
\bigl\|M^{[2]}_{t}F_{1,j}-F_{1,j}\bigr\|\leq
\biggl(\frac{17}{12}+\frac{27}{16\mathrm{e}}
\biggr)\frac{1}{t^{2}},\qquad
j=2,3,\ldots.
\end{equation}
Moreover, using Lemma
\ref{leboex}, we have
%
\begin{equation}\label{lemtm1}
\bigl\|M^{[2]}_{t}F_{1,1}-F_{1,1}\bigr\|\leq\biggl(\frac{1}{2}+\frac{1}{6
\mathrm{e}}+\frac{9}{4\mathrm{e}^{2}} \biggr)\frac{1}{t^{2}}\leq
\biggl(\frac{17}{12}+\frac{27}{16\mathrm{e}} \biggr)\frac{1}{t^{2}}.
\end{equation}
Now, let the general scale parameters be $a_{i}, i=1,\ldots,n$. We use
the fact that given $X$, a gamma random variable of scale parameter 1,
$X/a_{i}$ is a gamma random variable of scale parameter $a_{i}$, and,
therefore, using Proposition \ref{presca}, (\ref{lemtmj}) and
(\ref{lemtm1}), we have for each $a_{i}, i=1,\ldots, n$, and $j\in
\mathbb{N}^{*}$,
%
\begin{equation}\label{cotesc}
\bigl\|M_{t}^{[2]}F_{a_{i},j}-F_{a_{i},j}\bigr\|=\bigl\|
M_{t/a_{i}}^{[2]}F_{1,j}-F_{1,j}\bigr\|\leq
\biggl(\frac{17}{12}+\frac{27}{16\mathrm{e}} \biggr)\frac{a^{2}_{i}}{t^{2}}.
\end{equation}
Thus, using (\ref{linear}) and
(\ref{cotesc}), we have
\begin{eqnarray}\label{lemtmg}
\bigl\|M_{t}^{[2]}F-F\bigr\|&\leq&\sum_{i=1}^{n} \sum
_{j=1}^{\infty}p_{ij}\bigl\|M_{t}^{[2]}F_{a_{i},j}-F_{a_{i},j}\bigr\|
\nonumber\\[-8pt]\\[-8pt]
&\leq&
\biggl(\frac{17}{12}+\frac{27}{16e} \biggr)\frac{\sum_{i=1}^{n}
(\sum_{j=1}^{\infty}p_{ij})a_{i}^{2}}{t^{2}}.\nonumber
\end{eqnarray}
This completes the proof of Proposition \ref{procot}.
\end{pf}

As a consequence of the previous result, we can provide error bounds
for compound distributions (that is, distribution functions of random
sums, as in (\ref{ransum})) when the summands are mixtures of Erlang
distributions, as stated in the
following result.
\begin{corollary}\label{cocomp}
Let $G$ be the
distribution function of a random sum, as in \textup{(\ref{ransum})},
in which
the sequence of $(X_{i})_{i\in\mathbb{N}^{*}}$ has a common
distribution $\mathcal{ME}(a_{1},\ldots, a_{n}),~0<a_{1}<\cdots<a_{n}$, as defined in
\textup{(\ref{mixgav})}. Let $M_{t}^{[2]}$ be as in \textup{(\ref{a2gamo})}.
We have that
\[
\bigl\|M_{t}^{[2]}G-G\bigr\|\leq\biggl(\frac{17}{12}+\frac{27}{16\mathrm{e}} \biggr)\frac
{(1-G(0))a_{n}^2}{t^2}
.
\]
\end{corollary}

\begin{pf}
The proof is immediate, taking into account that
a mixture of Erlang distributions $\mathcal{ME}(a_{1},\ldots,
a_{n})$, $0<a_{1}<\cdots
<a_{n}$, can be expressed as in (\ref{mixgag}) and compound
distributions of these random variables are also mixtures of Erlang
distributions (see \cite{wiwoon}, page 106, with a slight modification
in the
coefficients, as we allow a point mass at 0), that is, we can write
\[
G(x)=q_{0}+\sum_{j=1}^{\infty}q_{j}F_{a_{n},j}(x),\qquad x\geq0,
\]
in which $\{q_{j}, j=0,1,\dots\}$ form a probability mass
function (obviously, $q_{0}=G(0)$). The result follows using the above
expression and Proposition \ref{procot}.
\end{pf}

The class of phase-type distributions, of great importance in
applied probability, can be expressed as mixtures of Erlang distributions.
A phase-type distribution is defined as the time until absorption in
a continuous-time Markov chain with one absorbent state (see, for
instance, \cite{larain}, Chapter II or \cite{asruin}, Chapter VIII, and the
references therein). A phase-type distribution can be expressed in
terms of a matrix exponential as follows. Consider a vector
$\mathbf{\alpha}=(\alpha_{1},\ldots,\alpha_{n})$ of non-negative
numbers such that $\alpha_{1}+\cdots+\alpha_{n}\leq1$. Let $A$ be a
$n\times n$ matrix with negative diagonal entries, non-negative
off-diagonal entries and non-positive row sums. A non-negative
random variable $X$ is a phase-type distribution
$\mathit{PH}(\mathbf{\alpha},A)$ if its distribution function can be written as
\[
F(x)=1-\mathbf{\alpha}\mathrm{e}^{xA}\mathbf{1}',\qquad x\geq0,
\]
in which $\mathbf{1}'$
represents the transpose of the $n$th dimensional vector
$\mathbf{1}=(1,\ldots,1)$. Note that phase-type distributions are
absolutely continuous random variables when
$\alpha_{1}+\cdots+\alpha_{n}=1$, having positive mass at 0 (of
magnitude $ 1-(\alpha_{1}+\cdots+\alpha_{n})$) when
$\alpha_{1}+\cdots+\alpha_{n}<1$. Phase-type distributions have been
extensively studied from both theoretical and practical points of
view. For instance, it is well known that phase-type distributions have rational
Laplace transforms, thus allowing numerical computation using our
approximation procedures. Also, in the next proposition, we will give
an expression of phase-type distributions in terms of mixtures of
Erlang distributions. This, together with Proposition \ref{procot},
provides our approximations with rates of convergence. The proof of the
next result is based on the following property of phase-type
distributions, due to Maier (see \cite{mathea}, page 591).
Let $f$ be the density of an absolutely continuous phase-type
distribution. There exists some $c>0$ verifying
%
\begin{equation}\label{phasede}
c_{j}:=\frac{\mathrm{d}^{j} }{\mathrm{d}x^{j}}\mathrm{e}^{cx}f(x)
\bigg|_{x=0}>0,\qquad
j\in\mathbb{N}.
\end{equation}
We are now in a position to
state the following.

\begin{proposition}\label{cophas}
Let $F$ be a phase-type distribution
$\mathit{PH}(\mathbf{\alpha},A)$, with $\alpha_{1}+\cdots+\alpha_{n}>0$. Let
$c>0$ be such that the absolutely continuous part of $F$ satisfies
the property
\textup{(\ref{phasede})}. Then, $F$ can
be expressed as a mixture of Erlang distributions, that is,
%
\begin{equation}\label{copha1}
F(x)=p_{0}+\sum_{j=1}^{\infty}p_{j}F_{c,j}(x),\qquad x\geq0,
\end{equation}
in which $p_{0}=1-(\alpha_{1}+\cdots+\alpha_{n})$.
\end{proposition}

\begin{pf}
To prove (a), assume first that $F$ is absolutely
continuous, that is, that $\alpha_{1}+\cdots+\alpha_{n}=1$.
Its density is then given by
$f(x)=-\mathbf{\alpha}\mathrm{e}^{xA}A\mathbf{1}', x>0 $. We choose a
$c>0$ verifying (\ref{phasede}). Note that we can write
%
\begin{equation}\label{maierp}
\mathrm{e}^{cx}f(x)=-\mathbf{\alpha}\mathrm{e}^{x(cI-A)}A\mathbf{1}',\qquad
x\geq0.
\end{equation}
It can be easily checked that
the function $-\mathbf{\alpha}\mathrm{e}^{x(cI-A)}A\mathbf{1}', x\in
\mathbb{R}$ is analytic in $\mathbb{R}$, so if we consider the Taylor
series expansion of this function around 0 and take into account
(\ref{phasede}) and (\ref{maierp}), we have
\[
\mathrm{e}^{cx}f(x)=\sum_{j=0}^{\infty}c_{j}\frac{x^{j}}{j!},\qquad x>0,
\]
from which we can write (recall (\ref{gamden}))
\[
f(x)=\sum_{j=0}^{\infty}\frac{c_{j}}{c^{j+1}}\frac
{c^{j+1}x^{j}\mathrm{e}^{-cx}}{j!}=\sum_{j=0}^{\infty}\frac
{c_{j}}{c^{j+1}}f_{c,j+1}(x),\qquad x>0,
\]
and, in this way, we obtain the expression of $f$ in terms of a
mixture of Erlang densities with shape parameter $c$ (by
construction, the coefficients are non-negative and integrating both
sides in the above expression, we see that their sum is 1). As a
consequence, we can write
%
\begin{equation}\label{acpart}
F(x)=\sum_{j=1}^{\infty}\frac
{c_{j-1}}{c^{j}}F_{c,j}(x),\qquad x\geq0,
\end{equation}
thus
having expressed $F$ as a mixture of Erlang distributions, as in (\ref
{copha1}).
Now, assume that $0<\alpha_{1}+\cdots+\alpha_{n}<1$. This means that
$F$ has a point mass at 0 of magnitude
$p_{0}:=1-(\alpha_{1}+\cdots+\alpha_{n})$. The absolutely continuous
part of $F$ ($F^{\mathrm{ac}}$) is a phase-type distribution
($\mathit{PH}(\bar{\mathbf{\alpha}},A)$) with
$\bar{\mathbf{\alpha}}=(\alpha_{1}+\cdots+\alpha_{n})^{-1}\mathbf
{\alpha}$.
Let $c>0$ be such that $F^{\mathrm{ac}}$ verifies property (\ref{maierp}). We can
write, thanks to (\ref{acpart}),
\[
F(x)=p_{0}+(1-p_{0})F^{\mathrm{ac}}(x)=p_{0}+\sum_{j=1}^{\infty}(1-p_{0})\frac
{c_{j-1}}{c^{j}}F_{c,j}(x),\qquad
x\geq0 .
\]
This completes the proof of Proposition
\ref{cophas}.
\end{pf}

\begin{remark}
Expansions similar to those given in Proposition \ref{cophas}
can be found in \cite{larain}, page~58. These expansions are obtained
using a representation $\mathit{PH}(\mathbf{\alpha},A)$ of the distribution
under consideration. Note that if we denote by $\|A\|$ the maximum
absolute value of the entries of $A$, then it is easy to check using
(\ref{maierp}) (see \cite{ocphas}, page 751) that $c=\|A\|$ verifies
(\ref{phasede}). However, as the representation of a phase-type
distribution is
not unique, this value might not be the optimum one. Also, observe
that the error bound given in (\ref{cotmix}) indicates that we
should take $c$ to be as small as possible. This problem, then, is closely
connected to Conjecture 6 in \cite{ocphas}, concerning the minimum
$c$ satisfying (\ref{phasede}) and its relation with a phase-type
representation having $\|A\|$ as small as possible. To the best of
our knowledge, this conjecture remains unsolved.
\end{remark}

\section*{Acknowledgments}
I would like to thank Jos\'{e} Garrido for suggesting the final
applications to phase-type distributions when I was at Concordia
University and to two anonymous referees for helpful comments. This
research has been partially supported by research Grants 2006-CIE-05
(University of Zaragoza), MTM2007-63683
and PR 2007-0295 (Spanish Government), E64 (DGA) and by FEDER
funds.\looseness=1

\printhistory


\begin{thebibliography}{00}

\bibitem{abwhla}
Abate, J. and Whitt, W. (1995). Numerical inversion of Laplace
transforms of probability distributions.
\textit{ORSA Journal of Computing} \textbf{7} 36--43.

\bibitem{abwhan}
Abate, J. and Whitt, W. (1996). An operational calculus for probability
distributions via Laplace Transforms.
\textit{Adv. in Appl. Probab.} \textbf{28} 75--113.
\MR{1372332}

\bibitem{abstha}
Abramowitz, M. and Stegun, I.A. (1964). \textit{Handbook of
Mathematical Functions}. Washington, DC: National Bureau of Standards.
\MR{0167642}

\bibitem{adcaon}
Adell, J.A. and de la Cal, J. (1993). On the uniform
convergence of normalized Poisson mixtures to their mixing
distribution. \textit{Statist. Probab. Lett.} \textbf{18} 227--232.
\MR{1241619}

\bibitem{adcaap}
Adell, J.A. and de la Cal, J. (1994). Approximating gamma
distributions by normalized negative binomial distributions. \textit{J.
Appl. Probab.} \textbf{31} 391--400.
\MR{1274795}

\bibitem{adsadi}
Adell, J.A. and Sang\"{u}esa, C. (1999).
Direct and converse inequalities for positive linear operators
on the positive semi-axis. \textit{J. Austral. Math. Soc. Ser. A}
\textbf{66} 90--103.
\MR{1658683}

\bibitem{alonso}
Alzer, H. (1997). On some inequalities
for the gamma and psi functions. \textit{Math. Comp.} \textbf{66}
373--389.
\MR{1388887}

\bibitem{asruin}
Asmussen, S. (2000). \textit{Ruin Probabilities}. Singapore: World
Scientific.
\MR{1794582}

\bibitem{emgrso}
Embrechts, P., Gr\"{u}bel, R. and Pitts, S.M. (1993). Some
applications of the fast Fourier transform algorithm in insurance
mathematics. \textit{Statist. Neerlandica} \textbf{47} 59--75.
\MR{1208036}

\bibitem{feanin}
Feller, W. (1971). \textit{An Introduction to Probability Theory
and Its Applications} \textbf{II}, 2nd ed. New York: Wiley.

\bibitem{grheco}
Gr\"{u}bel, R. and Hermesmeier, R. (2000). Computation of compound
distributions II: Discretization errors and Richardson
extrapolation. \textit{Astin Bull.} \textbf{30} 309--331.
\MR{1946567}

\bibitem{larain}
Latouche, G. and Ramaswami, V. (1999). \textit{Introduction to Matrix Analytic
Methods in Stochastic Modelling}. Philadelphia: ASA-SIAM.
\MR{1674122}

\bibitem{mathea}
Maier, R.S. (1991). The algebraic construction of phase-type
distributions. \textit{Comm. Statist. Stochastic Models} \textbf{7}
573--602.
\MR{1139070}

\bibitem{ocphas}
O'Cinneide, C.A. (1999). Phase-type distributions: Open
problems and a few properties. \textit{Comm. Statist. Stochastic
Models} \textbf{15} 731--757.
\MR{1708454}


\bibitem{qicuso}
Qi, F., Cui, R., Chen, C. and Guo, B. (2005). Some completely monotonic
functions involving
polygamma functions and an application. \textit{J. Math. Anal. Appl.}
\textbf{310} 303--308.
\MR{2160691}

\bibitem{saerro}
Sang\"{u}esa, C. (2008). Error bounds in approximations of
compound distributions using gamma-type operators. \textit{Insurance
Math. Econom.} \textbf{42} 484--491.
\MR{2404310}

\bibitem{saunif}
Sang\"{u}esa, C. (2008). Uniform error bounds in continuous
approximations of
nonnegative random variables using Laplace Transforms.
Preprint. Available at:
\url{http://www.unizar.es/galdeano/preprints/2008/prep08-01.pdf}.

\bibitem{surecu}
Sundt, B. (2002). Recursive evaluation of aggregate claims
distributions. \textit{Insurance Math. Econom.} \textbf{30} 297--322.
\MR{1921109}\

\bibitem{wiwoon}
Willmot, G.E. and Woo, J.K. (2007). On the class of Erlang
mixtures with risk theoretical applications. \textit{N. Am. Actuar. J.}
\textbf{11}
99--105.
\MR{2380721}

\end{thebibliography}
\end{document}